\documentclass[a4paper,12pt]{amsart}
\usepackage{amssymb}
\usepackage{ifthen}
\usepackage{cite}
 \usepackage[dvips]{graphicx}
\nonstopmode\numberwithin{equation}{section}
\usepackage{amssymb}
\usepackage{ifthen}
\usepackage{graphicx}
\usepackage{amsmath}
\usepackage[T1]{fontenc} %skandit
\usepackage[utf8]{inputenc}
\usepackage[usenames,dvipsnames]{color}
\usepackage{color}
\usepackage[english]{babel}
\usepackage{fancyhdr}
\usepackage{fancybox}
\usepackage{tikz}
\newtheorem*{theoA}{Theorem A}
\newtheorem*{theoB}{Theorem B}
\newtheorem*{theoC}{Theorem C}

\nonstopmode \numberwithin{equation}{section}
\theoremstyle{plain}

\newtheorem{conj}{Conjecture}

\theoremstyle{definition}
\newtheorem{defi}{Definition}[section]

\newtheorem{thm}{Theorem}[section]
\newtheorem{prob}{Problem}[section]
\newtheorem{cor}{Corollary}[section]

\newtheorem{prop}{Proposition}[section]
\newtheorem{rem}{Remark}[section]
\newtheorem{lem}{Lemma}[section]

\newtheorem{Open Ques}{Open question}[section]

\newcounter{minutes}
\divide\time by 60
\newcounter{hours}
\multiply\time by 60
\addtocounter{minutes}{-\time}
\newcounter {own}
\def\theown {\thesection       .\arabic{own}}

\newenvironment{pf}[1][]{%
 \vskip 3mm
 \noindent
 \ifthenelse{\equal{#1}{}}%
  {{\slshape Proof. }}%
  {{\slshape #1.} }%
 }%
{\qed\bigskip}

\newcounter{alphabet}

\def\be{\begin{equation}}
\def\ee{\end{equation}}

\newcommand{\bee}{\begin{enumerate}}
\newcommand{\eee}{\end{enumerate}}

\newcommand{\blem}{\begin{lem}}
\newcommand{\elem}{\end{lem}}
\newcommand{\bthm}{\begin{thm}}
\newcommand{\ethm}{\end{thm}}
\newcommand{\bcor}{\begin{cor}}
\newcommand{\ecor}{\end{cor}}
\newcommand{\beg}{\begin{examp}}
\newcommand{\eeg}{\end{examp}}
\newcommand{\begs}{\begin{examples}}
\newcommand{\eegs}{\end{examples}}
\usepackage{mathrsfs}
\newcommand{\bdefn}{\begin{defn}}
\newcommand{\edefn}{\end{defn}}

\newcommand{\bprob}{\begin{prob}}
\newcommand{\eprob}{\end{prob}}
\newcommand{\bei}{\begin{itemize}}
\newcommand{\eei}{\end{itemize}}

\newcommand{\bcon}{\begin{conj}}
\newcommand{\econ}{\end{conj}}
\newcommand{\bcons}{\begin{conjs}}
\newcommand{\econs}{\end{conjs}}
\newcommand{\bprop}{\begin{prop}}
\newcommand{\eprop}{\end{prop}}
\newcommand{\br}{\begin{rem}}
\newcommand{\er}{\end{rem}}
\newcommand{\brs}{\begin{rems}}
\newcommand{\ers}{\end{rems}}
\newcommand{\bo}{\begin{obser}}
\newcommand{\eo}{\end{obser}}
\newcommand{\bos}{\begin{obsers}}
\newcommand{\eos}{\end{obsers}}
\newcommand{\bpf}{\begin{pf}}
\newcommand{\epf}{\end{pf}} 
\newcommand{\ba}{\begin{array}}
\newcommand{\ea}{\end{array}}
\newcommand{\beq}{\begin{eqnarray}}
\newcommand{\beqq}{\begin{eqnarray*}}
\newcommand{\eeq}{\end{eqnarray}}
\newcommand{\eeqq}{\end{eqnarray*}}

\begin{document}

\title{Bohr's theorem for Ces\'aro operator and certain integral transforms over octonions}

\author{Molla Basir Ahamed$^*$}
\address{Molla Basir Ahamed, Department of Mathematics, Jadavpur University, Kolkata-700032, West Bengal,India.}
\email{mbahamed.math@jadavpuruniversity.in}

\author{Sabir Ahammed}
\address{Sabir Ahammed, Department of Mathematics, Jadavpur University, Kolkata-700032, West Bengal,India.}
\email{sabira.math.rs@jadavpuruniversity.in}

\subjclass[{AMS} Subject Classification:]{30A10, 30B10, 30G35, 30H05, 44A10, 44A55.}

\def\thefootnote{}
\footnotetext{ {\tiny File:~\jobname.tex,
printed: \number\year-\number\month-\number\day,
          \thehours.\ifnum\theminutes<10{0}\fi\theminutes }
} \makeatletter\def\thefootnote{\@arabic\c@footnote}\makeatother
\maketitle
\pagestyle{myheadings}
\markboth{M. B. Ahamed and S. Ahammed}{Bohr's theorem for Ces\'aro operator and certain integral transforms over octonions}
\begin{abstract}
	In this paper, we first establish the Bohr's theorem for Ces\'aro operator defined for  $f\in \mathcal{SRB}(\mathbb{B})$ of slice regular functions in the open unit ball $\mathbb{B}$ of the largest alternative division algebras of octonions $\mathbb{O}$, such that $|f(x)| \leq 1$ for all $x \in \mathbb{B}$. Next, we establish Bohr type inequalities for Bernardi operator for the functions $f\in \mathcal{SRB}(\mathbb{B})$, and with the help of this, we obtain Bohr type inequality for Libera operator and Alexander operator.	Finally, we obtain Bohr-type inequalities for certain integral transforms, namely Fourier (discrete) and Laplace (discrete) transforms for $f\in \mathcal{SRB}(\mathbb{B}).$ All the results are proven to be sharp.
\end{abstract}
\section{\bf Introduction}
In the early twentieth century, the study of Dirichlet series—series of the form $\sum_{k=0}^{\infty}a_kk^{-s}$, where $a_k$ and $s$ are the complex numbers—became a highly attractive research topic. Harald Bohr initiated investigations into their convergence, focusing primarily on determining the largest region in which a given Dirichlet series converges absolutely. The absolute convergence problem of Dirichlet series (see \cite[p. xxiv]{Defant-2019}) aims to determine the maximum possible width of the strip where the series converges uniformly but not absolutely. Equivalently, this problem can be reformulated as finding the supremum of the abscissa of absolute convergence for a Dirichlet series that converges on $\{z\in\mathbb{C} : {\rm Re}(z)>0\}$ with a bounded limit function.\vspace{2mm}

 Let $\mathcal{H}(\mathbb{D})$ denote the class of analytic functions $f(z)=\sum_{k=0}^{\infty}a_kz^k$  in the open unit disk $\mathbb{D}:=\{ z\in \mathbb{C}:|z|<1\}$ such that $|f(z)|\leq 1$ for all $z\in \mathbb{D}.$ In $1914$, Harald Bohr \cite{Bohr-1914} proved the following remarkable result concerning the universal constant $|z|=1/3$ for functions in $\mathcal{H}(\mathbb{D})$, which is now known as Bohr's theorem.
\begin{theoA} If $f(z)=\sum_{k=0}^{\infty}a_kz^k\in \mathcal{H}(\mathbb{D}),$  then 
	\begin{align}\label{Eq-1.11}
		\sum_{k=0}^{\infty}|a_k||z|^k\leq 1 \;\;\;\mbox{for}\;\;\; |z|\leq \dfrac{1}{3}.
	\end{align}
	The constant $1/3$ is best possible.
\end{theoA}
The inequality $\eqref{Eq-1.11}$ fails when $ |z|>{1}/{3} $ in the sense that there are functions in $ \mathcal{H}\left(\mathbb{D}\right) $ for which the inequality is reversed when $ |z|>{1}/{3} $. H. Bohr initially showed in \cite{Bohr-1914} the inequality \eqref{Eq-1.11} holds only for $|z|\leq1/6$, which was later improved independently by M. Riesz, I. Schur, F. Wiener and some others. The sharp constant $1/3$  and the inequality \eqref{Eq-1.11} in Theorem A are called respectively, the Bohr radius and the  classical Bohr inequality for the family $ \mathcal{H}\left(\mathbb{D}\right) $.  Several other proofs of this interesting inequality were given in different articles  (see \cite{Sidon-1927,Paulsen-Popescu-Singh-PLMS-2002,Tomic-1962}). \vspace{1.2mm}

In $1995$, Dixon \cite{Dixon-BLMS-1995} applied the Bohr inequality to the long-standing problem of characterizing Banach algebras that satisfy the von Neumann inequality. This application garnered considerable attention from various researchers, leading to investigations of the phenomenon in a variety of function spaces.However, it is important to note that not every class of functions exhibits the Bohr phenomenon. For instance, Bénéteau et al. \cite{Beneteau- Dahlner-Khavinson-CMFT-2004} showed that the Bohr phenomenon does not occur in the Hardy space $ H^p(\mathbb{D},X),$ where $p\in [1,\infty),$ and Hamada \cite{Hamada-Math.Nachr-2021} demonstrated that there is no Bohr radius for holomorphic mappings with values in the unit ball of a complex Hilbert space ${H}$ when $\dim {H}\geq 2.$Moreover, $\mathcal{H}\left(\mathbb{D}\right)$ is not the sole class of analytic functions for which Bohr radii are studied; various other function classes and certain integral operators have also been investigated. Some of these are discussed below.Recently, there have been several exciting advancements in this area. The study of the Bohr phenomenon for holomorphic functions in several complex variables and in Banach spaces has become particularly active (see \cite{Aizenberg-PAMS-2000,Aizenberg-Aytuna-Djakov-JMAA-2001, Galicer-Mansilla-Muro-TAMS-2020, Kumar-PAMS-2022, Lata-Singh-PAMS-2022, Liu-Ponnusamy-PAMS-2021, Hamada-RM-2024,Boas-Khavinson-PAMS-1997}). Research has also extended to the non-commutative setting (see e.g., \cite{Paulsen-Singh-BLMS-2006,Popescu-IEOT-2019}), including operator-valued generalizations of the Bohr theorem in the single-variable case \cite{Paulsen-Popescu-Singh-PLMS-2002}, with further references available in \cite{Popescu-TAMS-2007}.\vspace{1.2mm}

As a generalization of holomorphic functions of one complex variable, the theory of slice regular functions of one quaternionic variable was initiated by Gentili and Struppa \cite{Gentili-Struppa-AVM-2007} and further developed for Clifford algebras \cite{Colombo-Sabadini-Struppa-IJM-2009} and octonions \cite{Gentili-Struppa-RMJM-2010}. Based on the concept of stem functions, these three classes of functions were eventually unified and generalized into real alternative algebras \cite{Ghiloni-Perotti-AdvM-2011}. Following the historical path,   Rocchetta \emph{et al.} \cite{Rocchetta-Gentili-Sarfatti-MN-2012} have established Bohr's theorem for a special case of non-commutative (but associative) algebra of quaternions, where regular quaternionic rational transformations from \cite{Bisi-Gentili-IUMJ-2009,Stoppato-AG} are heavily utilized. This is because quaternions form a skew but associative field. When restricted to each slice, the function $ f$ described in Theorem B can be viewed, according to \cite[Lemma 3.2]{Xu-PRSE-2021}, as a vector-valued holomorphic function mapping from \( \mathbb{D} \) into the unit ball of \( \mathbb{C}^4 \) endowed with the standard Euclidean norm. Blasco \cite[Theorem 1.2]{Blasco-OTAA-2010} introduced the Bohr radius for holomorphic functions mapping from $\mathbb{D}$ into the unit ball of $\mathbb{C}^n$ (with $n \geq 2$) and proved that it is zero. From this perspective, the theory of slice regular functions differs from that of vector-valued holomorphic functions.\vspace{1.2mm} 

Recently, Xu \cite{Xu-PRSE-2021} established the Bohr's theorem  for the class $\mathcal{SRB}(\mathbb{B})$ of slice regular functions in the open unit ball $\mathbb{B}$ of the largest alternative division algebras of octonions $\mathbb{O}$, such that for all $x \in \mathbb{B}$, the condition $|f(x)| \leq 1$ is satisfied  and it's sharpness.
\begin{theoB}\emph{\cite[Theorem 1.3]{Xu-PRSE-2021}}
	Let $f(x)=\sum_{k=0}^{\infty}x^ka_k\in \mathcal{SRB}(\mathbb{B})$ with $a_k\in \mathbb{O}.$ Then 
	\begin{align*}
		\sum_{k=0}^{\infty}|x^ka_k|\leq 1\;\mbox{for}\; |x|\leq \dfrac{1}{3}.
	\end{align*}
	The constant $1/3$ is sharp.
\end{theoB}
 In $1931,$ Hardy and Littlehood \cite{Hardy-Littlewood-MZ-1931} studied the classical Ces\'aro operator (see for more information \cite{Siskakis-JLMS-1987,Siskakis-PAMS-1990}) and defined as
\begin{align}\label{Eq-1.1111}
	T[f](z):=\int_{0}^{1}\dfrac{f(tz)}{1-tz}dt =\sum_{k=0}^{\infty}\left(\dfrac{1}{k+1}\sum_{n=0}^{k}a_n\right)z^k,
\end{align}
where $f(z)=\sum_{k=0}^{\infty}a_kz^k$ is analytic in $\mathbb{D}$ and later, several authors have studied the boundedness of this operator on various function spaces (see \cite{Albanese- Bonet-Ricker-SM-2018}). A generalized form of the Ces\'aro operators and its Boundedness is studied in \cite{Agrawal-Howlett-Lucas-NaikPonnusamy-JCAM-2005}. Udez \emph{et al.} \cite{udez-Bonilla-Muller- Peris-JAM-2020}, extensively studied the Ces\'aro mean and boundedness of Ces\'aro operators  on Banach spaces and Hilbert spaces. For details of  Ces\'aro operators and related properties, we refer to the article\cite{Stempak-PRSE-1994}. Kayumov \emph{et al.} \cite{Kayumov-Khammatova-Ponnusamy-CRMA-2020,Kayumov-Khammatova-Ponnusamy-MJM-2020} established the Bohr's theorem for Ces\'aro operator defined for analytic function in $\mathbb{D}$ and  established Bohr-Rogosinski phenomenon for ces\'aro operators in \cite{Kayumov-Khammatova-Ponnusamy-JMAA-2021}. Kumar and Sahoo \cite{Kumar-Sahoo-MJM-2021} studied Bohr-type inequalities for more generalized integral operators. In \cite{Allu-Halder-PEMS-2023, Aha-Aha-MJM-2025}, the authors have studied Bohr-type inequality and Bohr-Rogosinski inequality for Ces\'aro operator defined for operator-valued holomorphic functions in a simply connected domain containing the unit disk $\mathbb{D}.$ As noted in \cite{Kayumov-Khammatova-Ponnusamy-CRMA-2020},
\[ \bigg|T[f](z)\bigg|\leq \dfrac{1}{r}\log\dfrac{1}{1-r}\]
for each $|z|=r<1.$ On the other hand, from \eqref{Eq-1.1111}, we have the obvious estimate 
\[ \bigg|T[f](z)\bigg|\leq \sum_{k=0}^{\infty}\left(\dfrac{1}{k+1}\sum_{n=0}^{k}|a_n|\right)|z|^k, \]
the absolute sum of the series \eqref{Eq-1.1111}.  Kayumov \emph{et al.} \cite{Kayumov-Khammatova-Ponnusamy-CRMA-2020} determined the precise radius $r$ for which the absolute sum shares the same upper bound $(1/r) \log(1/(1-r))$. This investigation was significant because a convergent series is not necessarily absolutely convergent. In fact, they established the following sharp inequality.
\begin{theoC}\cite[Theorem 1]{Kayumov-Khammatova-Ponnusamy-CRMA-2020}\label{Thm-1.15}
	If $f(z)=\sum_{k=0}^{\infty}a_kz^k\in \mathcal{H}(\mathbb{D}),$ then 
	\[ \sum_{k=0}^{\infty}\left(\dfrac{1}{k+1}\sum_{n=0}^{k}|a_n|\right)r^k\leq \dfrac{1}{r}\log\dfrac{1}{1-r}\]
	for $|z|=r\leq R\approx0.5335,$ where $R$ is the unique root in $(0,1)$ of the equation 
	\begin{align}\label{Eq-1.12}
		2y-3(1-y)\log\dfrac{1}{1-y}=0
	\end{align}
	that cannot be improved. 
\end{theoC}
The Bohr radius problem has been a subject of extensive research in recent years, focusing primarily on classes of self-analytic maps on $\mathbb{D}$, harmonic mappings, and holomorphic functions of several complex variables. While the literature includes numerous studies on this problem, such as those concerning improved Bohr inequalities and their refinement, there remains a significant gap: the corresponding study involving the Cesàro operator and certain integral transforms has not garnered comparable attention, with only a few results published (see \cite{Aha-Aha-MJM-2025,Kayumov-Khammatova-Ponnusamy-JMAA-2021,Kayumov-Khammatova-Ponnusamy-MJM-2020,Kumar-Sahoo-MJM-2021,Kumar-CVEE-2023,Ong-N-IJS-2024}). This gap is the central motivation for the current paper. We aim to contribute to the understanding of the Bohr radius for the Cesàro operator and certain integral transforms defined on the space over octonions and establish the sharpness of the resulting inequalities.\vspace{1.2mm}

Inspired by the work of \cite{Ghiloni-Perotti-AdvM-2011}, we pose the following natural question to further investigate Bohr's theorem. 
\begin{prob}\label{Q-1.1}
	Can we establish Bohr-type inequalities for Ces\'aro operators, Fourier (discrete) and Laplace (discrete) transforms defined on the space of slice regular functions over octonions?
\end{prob}
In this paper, our objective is to give affirmative answer to Problem \ref{Q-1.1} and establish results. The organization of the paper is follows: In Section \ref{Sec-2}, we state the main results and some remarks. In Section \ref{Sec-3}, we recall necessary definitions and preliminary results used in the sequel for slice regular functions. In Section \ref{Sec-4}, we give the details of the proof of the main results.
\section{\bf Main results}\label{Sec-2}
For  $f(z)=\sum_{k=m}^{\infty}a_kz^k$ is complex-valued analytic function in $\mathbb{D}$, Miller and Mocanu \cite[P. 11]{Miller-Mocanu-2000}, defined the Bernardi operator as follows: 
\begin{align*}
	L_\gamma[g](z):=\sum_{k=m}^{\infty}\dfrac{a_k}{k+\gamma}z^k=\int_{0}^{1}g(tz)t^{\gamma-1}dt,
\end{align*}
where $\gamma>-m$ and $m\geq0$ is an integer. The function $L_\gamma[g]$ is analytic in $\mathbb{D}$ and several properties of $L_\gamma[g]$  when $m=1$ with normalization, are well-known (see \cite{Miller-Mocanu-2000,Ponnusamy-Sahoo-JMAA-2008}). For $\gamma=1$ and $m=0,$ we obtain the well-known Libera operator (see \cite{Miller-Mocanu-2000}) defined by 
\begin{align*}
	L[g](z):=2\sum_{k=0}^{\infty}\dfrac{a_k}{k+1}z^k=2\int_{0}^{1}g(tz)dt, \;\;\mbox{for}\;\; z\in \mathbb{D}.
\end{align*}
For $\gamma=0,$ $m=1$ and $f(z)=\sum_{k=1}^{\infty}a_kz^k,$ we obtain the well-known Alexander operator (see \cite{Duren-1983,Miller-Mocanu-2000}) defined by
\begin{align*}
	J[g](z):=\int_{0}^{1}\dfrac{g(tz)}{t}dt=\sum_{k=1}^{\infty}\dfrac{a_k}{k}z^k, \;\;\mbox{for}\;\; z\in \mathbb{D},
\end{align*}
which has been extensively studied in the univalent function theory.\vspace{1.2mm}

Let $\{z_k\}^{N-1}_{k=0}$ be a sequence of complex numbers. The discrete Fourier transform (see \cite{Bachman-Narici-2000}, Chapter 6) is defined as 
\[ \mathcal{F}(z_n):= \sum_{k=0}^{N-1}z_k  e^{-{2\pi ikn}/{N}}. \]
For $f(z)=\sum_{n=0}^{\infty}a_nz^n,$ the discrete Fourier transform on the
coefficients $a_k$ from $k=0$ to $n$, which gives
\[ \mathcal{F}[f](z):=\sum_{n=0}^{\infty}\left(\sum_{k=0}^{n}a_ke^{-\frac{2\pi ink}{(n+1)}}\right)z^n.\]
Consider the discrete Laplace transform for a sequence $\{y_n\}^\infty_{n=0}$ as follows:
\[\mathcal{L}(y_n):=(x_k)=\sum_{n=0}^{\infty}\dfrac{y_n}{(k+1)^{n+1}},\]
wherever the series on the right-hand side converges  from \cite{Ameen-Hasan-Jarad-RNA-2019,Holm-CMA-2011}.
To obtain a Bohr-type inequality in the unit disk $\mathbb{D}$, Ong and Ng \cite{Ong-N-IJS-2024} defined the majorant series for the discrete Fourier transform  and the discrete Laplace transform respectively as follows:
\[ \mathcal{F}_f(r):=\sum_{n=0}^{\infty}\left(\sum_{k=0}^{n} \bigg|a_ke^{-\frac{2\pi ink}{(n+1)}}\bigg|\right)r^n,\]

and 
\[\mathcal{L}_f(r):=\sum_{n=0}^{\infty}\left( \sum_{k=0}^{n}\dfrac{|a_k|}{(n+1)^{k+1}}\right)r^n, \]
where $|z|=r<1,$ and using this, Ong and Ng obtained Bohr-type inequalities in \cite{Ong-N-IJS-2024}. \vspace{1.2mm}

In the same spirit of the definitions, we define $\beta$-Ces\'aro operator $(\beta>0)$ on the space of slice regular functions  $f(x)=\sum_{k=0}^{\infty}x^ka_k\in \mathcal{SRB}(\mathbb{B})$ by
\[ T^*_\beta[f](x):= \sum_{k=0}^{\infty}x^k\left(\dfrac{1}{k+1}\sum_{n=0}^{k}\dfrac{\Gamma(k-n+\beta)}{\Gamma(k-n+1)\Gamma(\beta)}a_n\right)=\int_{0}^{1}(1-tx)^{-\bullet}f(tx)dt,\]
 and the Bernardi operator defined as 
\[ L^*_\gamma[g](x):=\sum_{k=m}^{\infty}x^k\dfrac{a_k}{k+\gamma}=\int_{0}^{1}g(tx)t^{\gamma-1}dt,\]
for slice regular functions $g(x)=\sum_{k=m}^{\infty}x^ka_k \in \mathcal{SRB}(\mathbb{B})$ and $\gamma>-m,$ here $m\geq0$ is an integer. With the help of the Bernardi operator we also obtain the Bohr radii for some known operators.\vspace{1.2mm}

In the following, to derive a sharp Bohr-type inequality, the discrete Fourier transform applied to slice regular functions over octonions, for $f(x)=\sum_{k=0}^{\infty}x^ka_k,$ $x\in \mathbb{B}_J,$ where $J\in \mathbb{S},$ we perform the discrete Fourier transform on the coefficients $a_k$ from $k=0$ to $n$, which gives
\[ \mathcal{F}^*[f](x):=\sum_{n=0}^{\infty}x^n\left(\sum_{k=0}^{n}a_ke^{-\frac{2\pi Jnk}{(n+1)}}\right)\] 
and the associate majorant series of $\mathcal{F}[f](x)$ as follows:
\[ \mathcal{F}^*_f(|x|):=\sum_{n=0}^{\infty}|x|^n\left(\sum_{k=0}^{n} \bigg|a_ke^{-\frac{2\pi Jnk}{(n+1)}}\bigg|\right),\]
where $|x|<1.$\vspace{1.2mm}

 Before answering Problem \ref{Q-1.1}, we need to establish some preliminaries. For $f\in \mathcal{SRB}(\mathbb{B})$ and $\beta>0,$ an elementary estimation of the absolute value of the integral yields a sharp inequality
\[ \bigg|T^*_\beta[f](x)\bigg|\leq \begin{cases}
	\dfrac{1}{r}\left(\dfrac{1-(1-r)^{1-\beta}}{1-\beta}\right)\;\;\mbox{if}\;\; \beta\neq 1, \vspace{2mm}\\
	\dfrac{1}{r}\log\dfrac{1}{1-r}\;\;\;\;\;\;\;\;\;\;\;\;\;\;\;\;\;\mbox{if}\;\; \beta= 1,
\end{cases}\]
for each $|x|=r<1.$ 	Similar to the Bohr-type inequality for the operator $ T^*_\beta$, $\beta>0$, we also study the  Bohr-type inequality for the Bernardi operator $L^*_\gamma$. It is easy to calculate the following sharp bound
\begin{align*}
	|L^*_\gamma[f](x)|\leq \dfrac{1}{m+\gamma}|x|^m,\;\;|x|<1
\end{align*}
for $f(x)=\sum_{k=m}^{\infty}x^ka_k.$\vspace{2mm}

Now, we establish the Bohr-type inequality for the Ces'{a}ro operator defined on the space of slice regular functions over octonions. Our result is presented below.
\begin{thm}\label{Thm-1.16}
	Let $f(x)=\sum_{k=0}^{\infty}x^ka_k\in \mathcal{SRB}(\mathbb{B})$ and $0<\beta\neq1.$ Then 
	\begin{align*}
		\sum_{k=0}^{\infty}|x|^k\left(\dfrac{1}{k+1}\sum_{n=0}^{k}\dfrac{\Gamma(k-n+\beta)}{\Gamma(k-n+1)\Gamma(\beta)}|a_n|\right)\leq \dfrac{1}{r}\left(\dfrac{1-(1-r)^{1-\beta}}{1-\beta}\right)
	\end{align*} 	
	for $|x|=r\leq R_\beta,$ where $R_\beta$ is the positive root of the equation 
	\begin{align}\label{Eq-1.13}
		\dfrac{3(1-(1-y)^{1-\beta})}{1-\beta}-\dfrac{2\left((1-y)^{-\beta}-1\right)}{\beta}=0.
	\end{align}
	The radius $R_\beta$ cannot be improved. 
\end{thm}
	By considering the limiting case as $\beta\to 1$ in Theorem \ref{Thm-1.16}, we arrive at the following corollary, which provides an analog of Theorem C.
	\begin{cor}
	Let $f(x)=\sum_{k=0}^{\infty}x^ka_k\in \mathcal{SRB}(\mathbb{B})$. Then 
	\[ \sum_{k=0}^{\infty}|x|^k\left(\dfrac{1}{k+1}\sum_{n=0}^{k}|a_n|\right)\leq \dfrac{1}{r}\log\dfrac{1}{1-r}\]
	for $|x|=r\leq R\approx0.5335,$ where $R$ is the unique root in $(0,1)$ of the equation \eqref{Eq-1.12} and	that cannot be improved. 
	\end{cor}
\begin{rem}
	 Let $h(x)=\sum_{k=1}^{\infty}x^ka_k\in \mathcal{SRB}(\mathbb{B})$ be a slice regular function normalized by $h(0)=0.$ Then we can write $h(x)=xg(x),$ where $g(x):=\sum_{k=0}^{\infty}x^ka_{k+1}\in \mathcal{SRB}(\mathbb{B}).$ We have
	\[ \mathcal{C}_\beta[h](x):=\sum_{k=0}^{\infty}x^{k+1}\left(\dfrac{1}{k+1}\sum_{n=0}^{k}\dfrac{\Gamma(k-n+\beta)}{\Gamma(k-n+1)\Gamma(\beta)}\right)=xT^*_\beta[g](x).\]
	By using Theorem \ref{Thm-1.16}, we obtain
	\begin{align*}
		\sum_{k=0}^{\infty}|x|^{k+1}\left(\dfrac{1}{k+1}\sum_{n=0}^{k}\dfrac{\Gamma(k-n+\beta)}{\Gamma(k-n+1)\Gamma(\beta)}|a_{k+1}|\right)\leq\dfrac{1-(1-r)^{1-\beta}}{1-\beta},\;\;0<\beta\neq 1
	\end{align*}
	for $r\leq R_\beta,$ where $R_\beta$ is the positive root of the equation \eqref{Eq-1.13} that cannot be improved. Moreover, in the limit $\beta\to 1,$ we can indeed obtain the Bohr radius problem: If $h(x)=\sum_{k=1}^{\infty}x^ka_k\in \mathcal{SRB}(\mathbb{B})$ be a slice regular function normalized by $h(0)=0,$ then 
	
	\[ \sum_{k=0}^{\infty}|x|^{k+1}\left(\dfrac{1}{k+1}\sum_{n=0}^{k}|a_{n+1}|\right)\leq \log\dfrac{1}{1-r}\]
	for $r\leq R\approx0.5335,$ where $R$ is the unique root in $(0,1)$ of the equation \eqref{Eq-1.12} that cannot be improved.
\end{rem}
Next, we establish the Bohr-type inequality for the Bernardi operator, which is presented below as a sharp result.
\begin{thm}\label{Thm-1.17}
	Let $\gamma>-m,$ $m\geq0$ is an integer and $f(x)=\sum_{k=m}^{\infty}x^ka_k\in \mathcal{SRB}(\mathbb{B}).$ Then 
	\[ \sum_{k=m}^{\infty}\dfrac{1}{k+\gamma}|x^ka_k|\leq \dfrac{1}{m+\gamma}|x|^m\]
	for $|x|\leq R_\gamma,$ where $R_\gamma$ is is the positive root of the equation 
	
	\begin{align}\label{EEq-1.15}
		\dfrac{y^m}{m+\gamma}-2\sum_{k=m+1}^{\infty}\dfrac{y^k}{k+\gamma}=0
	\end{align} and $R_\gamma$ cannot be improved.
\end{thm}
By setting $\gamma = 1$ and $m = 0$ in Theorem \ref{Thm-1.17}, we obtain the Bohr radius for the Libera operator and the Alexander operator, respectively, as a special case.
\begin{cor}
	If	$f(x)=\sum_{k=0}^{\infty}x^ka_k\in \mathcal{SRB}(\mathbb{B}),$ then
	\[ \sum_{k=0}^{\infty}\dfrac{|x^ka_k|}{k+1}\leq 1,\;\;\mbox{for}\;\; r\leq R_1,\]
	where $R_1=0.5828\dots$ is the root of the equation $3y+2\log (1-y)=0.$ The constant $R_1$ cannot be improved. 
\end{cor}
\begin{cor}
	If	$f(x)=\sum_{k=1}^{\infty}x^ka_k\in \mathcal{SRB}(\mathbb{B}),$ then
	\[ \sum_{k=1}^{\infty}\dfrac{|x^{k-1}a_k|}{k}\leq 1,\;\;\mbox{for}\;\; r\leq R_1,\]
	where $R_1=0.5828\dots$ is the root in $(0, 1)$ of the equation $3y+2\log (1-y)=0.$ The constant $R_1$ cannot be improved. 
\end{cor}
In the next result, we establish a Bohr-type inequality for the discrete Fourier transform for $f\in \mathcal{SRB}(\mathbb{B})$.
\begin{thm}\label{Thm-1.19}
	If $f(x)=\sum_{k=0}^{\infty}x^ka_k\in \mathcal{SRB}(\mathbb{B}),$ then
\begin{align*}
	 \mathcal{F}^*_f(|x|)\leq \dfrac{1}{1-|x|}\;\;\;\mbox{for}\;\; |x|\leq \frac{1}{3}.
\end{align*}
	The constant $1/3$ cannot be improved.
\end{thm}
In the following result, we establish a Bohr-type inequality for the discrete Laplace transform.
\begin{thm}\label{Thm-1.20}
	If $f(x)=\sum_{n=0}^{\infty}x^na_n\in \mathcal{SRB}(\mathbb{B}),$ then 
\begin{align}\label{Eq-3.1}
	\mathcal{L}^*_f(|x|):=\sum_{n=0}^{\infty}|x|^n\left(\sum_{k=0}^{n}\dfrac{|a_k|}{(n+1)^{k+1}}\right)\leq \frac{1}{|x|}\log\dfrac{1}{1-|x|}
\end{align}
	for all $0<|x|< 1.$
\end{thm}

\subsection{A remark on Theorem \ref{Thm-1.20}}
The following estimate
\begin{align}\label{Eq-3.2}
	-\sum_{n=1}^{\infty}\frac{|x|^n}{n(n+1)^{n+1}}<0,
\end{align} is used in the proof of Theorem \ref{Thm-1.20}. In fact, the proof can be modified without relying on the inequality \eqref{Eq-3.2}. For instance, by utilizing the given estimate, we can proceed as follows:
\begin{align*}
	\sum_{k=1}^{n}\frac{1}{(n+1)^{k+1}}\leq \sum_{k=1}^{n}\frac{1}{(n+1)^2}=\frac{n}{(n+1)^2}.
\end{align*}
Thus, we have
\begin{align*}
	\sum_{n=0}^{\infty}\left(\sum_{k=1}^{n}\frac{1}{(n+1)^{k+1}}\right)|x|^n&\leq \sum_{n=0}^{\infty}\left(\sum_{k=1}^{n}\frac{n}{(n+1)^{2}}\right)|x|^n\\&=\frac{1}{|x|}\ln\left(\frac{1}{1-|x|}\right)-\frac{{\rm Li}_2(|x|)}{|x|}.
\end{align*}
Setting $|a_0|=p$ and $|x|=r,$  we obtain
\begin{align}\label{Eq-3.3}
	\mathcal{L}^{*}_f(r)\leq \frac{p}{r}\ln\left(\frac{1}{1-r}\right)+(1-p^2)\left(\frac{1}{r}\ln\left(\frac{1}{1-r}\right)-\frac{{\rm Li}_2(r)}{r}\right):=\mathcal{W}_{r}(p).
\end{align}
A simple computation gives us
\begin{align*}
	\mathcal{W}^{\prime}_{r}(p)=\frac{1}{r}\ln\left(\frac{1}{1-r}\right)-2p\left(\frac{1}{r}\ln\left(\frac{1}{1-r}\right)-\frac{{\rm Li}_2(r)}{r}\right)
\end{align*}
and 
\begin{align*}
	\mathcal{W}^{\prime\prime}_{r}(p)=-2 \left(\frac{1}{r}\ln\left(\frac{1}{1-r}\right)-\frac{{\rm Li}_2(r)}{r}\right)\leq 0
\end{align*}
for all $r\in (0, 1)$.\vspace{2mm}

\noindent Hence, we see that $\mathcal{W}^{\prime}_{r}$ is a decreasing function of $p$ in $[0, 1]$. Consequently, we have
\begin{align*}
	\mathcal{W}^{\prime}_{r}(p)\geq 	\mathcal{W}^{\prime}_{r}(1))=\frac{1}{r}\ln\left(\frac{1}{1-r}\right)-2\left(\frac{1}{r}\ln\left(\frac{1}{1-r}\right)-\frac{{\rm Li}_2(r)}{r}\right)\geq 0
\end{align*}
for all $p\in [0, 1]$. This shows that $\mathcal{W}_{r}$ is an increasing function of $p$ in  $[0, 1]$. As a result, we obtain 
\begin{align*}
	\mathcal{L}^{*}_f(r)\leq \mathcal{W}_{r}(p)\leq \mathcal{W}_{r}(1)=\frac{1}{r}\ln\left(\frac{1}{1-r}\right)\; \mbox{for all}\; r\in (0, 1).
\end{align*}
Our observation is that regardless of how we arrange the series, the upper bound of $\mathcal{L}^{*}_f(r)$ will remain the same due to the presence of the factor $(1-p^2)$ in the expression of the function $\mathcal{W}_{r}(p)$, as we take $p\to 1^{-}$.\vspace{2mm}

In contrast, our study shows that it is possible to find it by rearranging the terms using a different approach.
\begin{cor}\label{Cor-2.3}
	Under the assumption of Theorem \ref{Thm-1.20}, the inequality \eqref{Eq-3.1} holds for $r\leq r^*<1$, where $r^*$ is a root in $(0, 1)$ of the equation $\mathcal{P}(r)=0$, where
	\begin{align*}
		\mathcal{P}(r):=\frac{1}{r}\ln(1-r)+2\left(\frac{1}{r}\ln\left(\frac{1}{1-r}\right)-\frac{{\rm Li}_2(r)}{r}\right).
	\end{align*}
	 The inequality \eqref{Eq-3.1} is sharp.
\end{cor} 
\section{\bf Preliminaries}\label{Sec-3}
In this section, we recall necessary definitions and preliminary results used in the sequel for slice regular functions from \cite{Ghiloni-Perotti-AdvM-2011}.
\subsection{The algebra of octonions}
Let $\mathbb{C},$ $\mathbb{H},$ $\mathbb{O}$ denote the algebra of complex numbers, quaternions and octonions, respectively and$\{ 1, i,j,k\}$ be the standard basis of the non-commutative, associative, real algebra of quaternions with the multiplication rules 
\begin{align*}
	i^2=j^2=k^2=ijk=-1.
\end{align*} 
For each element $a=x_0+x_1i+x_2j+x_3k$ in  $\mathbb{H}$ $(x_0,x_1,x_2,x_3\in \mathbb{R}),$ the conjugate of $a$ is defined as   $\overline{a}=x_0-x_1i-x_2j-x_3k.$ By the well-known Cayley-Dickson process, the real algebra of octonions can be built from $\mathbb{H}$ as $\mathbb{O}=\mathbb{H}+l\mathbb{H}$ with 
\begin{align*}
	&\overline{a+lb}=\overline{a}-lb, \;\; (a+lb)+(c+ld)=(a+c)+l(b+d),\;\;\\& (a+lb)(c+ld)=(ac-d\overline{b})+l(\overline{a}d+cb),\;\;\mbox{for}\;\;a,b,c,d\in \mathbb{H}.
\end{align*} 
As a consequence, $\{ 1,i,j,k,li,lj,lk\}$ form the canonical real vector basis of $\mathbb{O}$. Every element $x\in \mathbb{O}$ can be composed into the real part $\mathrm{Re}(x)=(x+\overline{x})/2$ and the imaginary part $\mathrm{Im}(x)=x-\mathrm{Re}(x).$ Define the modulus of $x$ $|x|=\sqrt{x\overline{x}},$ which is exactly the usual Euclidean norm in $\mathbb{R}^8.$ Also, the modulus is multiplicative, i.e., $|xy|=|x||y|$ for all $x,y\in \mathbb{O}.$ Every non-zero element $x\in \mathbb{O}$ has a multiplicative inverse given by $x^{-1}=|x|^{-2}\overline{x}.$ The construction above shows that $\mathbb{O}$ is a non-commutative, non-associative, normed and division algebra. See for instance \cite{Baez-BAMS-2002} for more explanation on the octonions.\vspace{1.2mm} 

The set of square roots of $-1$ in $\mathbb{O}$ is the six-dimensional unit sphere given by 
\[ \mathbb{S}:=\{q\in \mathbb{H}: q^2=-1 \}.\] 

For each $I\in \mathbb{S},$ denote by $\mathbb{C}_I:=<1,I>\cong \mathbb{C}$ the subalgebra $\mathbb{O}$ of  generated over $\mathbb{R}$ by $1$ and $I$.\vspace{1.2mm}

Notice that each $x\in \mathbb{O}$ can be expressed as $x = \alpha + \beta I_x$ with $\alpha \in \mathbb{R},$ $\beta\in \mathbb{R}^+$ and $I_x \in \mathbb{S}$. This inconspicuous observation allows decomposing $\mathbb{O}$ into `complex slices'
\[ \mathbb{O}=\bigcup_{I\in \mathbb{S}}\mathbb{C}_I,\]
which derives the remarkable notion of slice regularity over octonions.

\subsection{Slice functions}
Given an open set ${D}$ of $\mathbb{C}$, invariant under the complex conjugation, its circularization $\Omega_D$ is defined by
\[ \Omega_D=\bigcup_{I\in \mathbb{S}}\{ \alpha+\beta I: \exists\; \alpha, \beta \in \mathbb{R},\; s.t.\; z= \alpha+i \beta\in D\}.\]
A subset $\Omega$ in $\mathbb{O}$ is called to be circular if $\Omega = \Omega_D$ for some $D\subset \mathbb{C}$. The open
unit ball $ \mathbb{B}=\{ x\in \mathbb{O}:|x|<1\}$ and the right half-space $\{x\in \mathbb{O}: \mathrm{Re}(x)>0 \}$ are two
typical examples of the circular domain. For a series expansion for slice regular functions and related properties, we refer to the article \cite{Stoppato-AM-2012}.

\begin{defi}
	A function $F:D\to \mathbb{O}\otimes_\mathbb{R}\mathbb{C}$ on an open set $D\subset \mathbb{C}$ invariant
	under the complex conjugation is called a stem function if the $\mathbb{O}$-valued components $F_1,$ $F_2$ of $F=F_1+iF_2$ satisfies
	\[ F_1(\overline{z})=F_1(z),\;\; F_2(\overline{z})=-F_2(z),\;\forall\; z=\alpha+i\beta.\]
\end{defi}
Each stem function $F$ induces a (left) slice function $f = \mathcal{I}(F):\Omega_D\to \mathbb{O}$ given by
\[ f(x):=F_1(z)+IF_2(z),\;\forall\; x=\alpha+I\beta\in \Omega_D.\]
We will denote the set of all such induced slice functions on $\Omega_D$ by
\[ \mathcal{S}(\Omega_D):=\{ f=\mathcal{I}(F): F
\; \mbox{is}\; \mbox{a}\; \mbox{stem}\; \mbox{function}\; \mbox{on}\; D\}.\]
Each slice function $f$ is induced by a unique stem function $F$ since $F_1$ and $F_2$ are determined by $f$. In fact, it holds that \[ F_1(z)=\dfrac{1}{2}\left(f(x)+f(\overline{x})\right),\;\; z\in \Omega_D,\]
and \[ f(z)=\begin{cases}
	\dfrac{1}{2I_x}\left( f(x)-f(\overline{x})\right)\;\;\mbox{if}\;\;z\in \Omega_D\setminus \mathbb{R},\\\vspace{4mm}
	0,\;\;\;\;\;\;\;\;\;\;\;\;\;\;\;\;\;\;\;\;\;\;\;\;\;\;\; \mbox{if}\;\; z\in \Omega_D\cap \mathbb{R}.
\end{cases}\]
Recall that a $\mathrm{C}^1$ function $F:D\to \mathbb{O}\otimes_\mathbb{R}\mathbb{C}$  is holomorphic if, and only if, its components $F_1$, $F_2$ satisfy the Cauchy-Riemann equations 
\[ \dfrac{\partial F_1}{\partial\alpha}=\dfrac{\partial F_2}{\partial\beta},\;\;\dfrac{\partial F_1}{\partial\beta}=-\dfrac{\partial F_2}{\partial\alpha},\;\; z=\alpha+i\beta\in D. \]
\begin{defi}
	A (left) slice function $f=\mathcal{I}(F)$ on $\Omega_D$ is regular if its stem function $F$ is holomorphic on $D.$ Denote the class of slice regular functions on $\Omega_D$ by
	\[ \mathcal{SR}\left(\Omega_D\right):=\{f=\mathcal{I}(F)\in \mathcal{S}(\Omega_D): F\;\mbox{is}\;\mbox{holomorphic}\;\mbox{on}\; D\}.\]
\end{defi}
For $f\in \mathcal{SR}\left(\Omega_D\right),$ the slice derivative is defined to be the slice regular function $f^\prime$
on $\Omega_D$ obtained as
\[ f^\prime(x):= \mathcal{I}\left( \dfrac{\partial F}{\partial z}(z)\right)=\dfrac{1}{2}\mathcal{I}\left(\dfrac{\partial F}{\partial \alpha}(z)-i\dfrac{\partial F}{\partial \beta}(z) \right). \]
Recall that $\mathbb{O}$ is non-associative but alternative, \emph{i.e.}, the associator $(x,y,z):= (xy)z-x(yz)$ of three elements $x,y,z\in \mathbb{O}$ is an alternating function in its arguments. Meanwhile, the Artin theorem asserts that the subalgebra generated by two elements of $\mathbb{O}$ is associative. Hence, a class of examples of slice regular functions is given by polynomials of one octonionic variable with coefficients in $\mathbb{O}$ on the right side. Indeed, each slice regular function $f$ defined in $\mathbb{B}$ admits the expansion of convergent power series 
\[ f(x)=\sum_{k=0}^{\infty} x^ka_k,\;\;\{ a_k\}\subset \mathbb{O},\] for all $x\in \mathbb{B}.$\vspace{1.2mm}

For simplicity, let $\mathbb{B}_I$ the intersection $\mathbb{B}\cap\mathbb{C}_I$ for any $I\in \mathbb{S}$. Then the restriction $f|_{\mathbb{B}_I}$  is holomorphic on $\mathbb{B}_I$. Furthermore, the relation between slice regularity and complex holomorphy can be presented as follows.
\begin{lem}(Splitting lemma)
	Let $\{ I_0=1,I,I_1,II_1,I_2, II_2,I_3,II_3\}$ be a splitting
	basis for $\mathbb{O}.$ For $f\in \mathcal{SR}(\Omega_D),$ there exist holomorphic functions $f_n:\Omega_D\cap\mathbb{C}_I\to \mathbb{C}_I,$ $m\in \{0,1,2,3\},$ such that 
	\[ f(z)=\sum_{n=0}^{3}f_n(z)I_n,\;\; \forall\; z\in \Omega_D\cap\mathbb{C}_I.\] 
\end{lem}
Due to that the pointwise product of two slice functions is not a slice function generally, the notion of slice product was introduced.
\begin{defi}
	Let $f=\mathcal{I}(F)$ and $g=\mathcal{I}(G)$ be in $\mathcal{S}(\Omega_D)$ with stem functions $F=F_1+iF_2$ and $G=G_1+iG_2.$ Then $FG=F_1G_1-F_2G_2+i\left(F_1G_2+F_2G_1\right)$ is
	still a stem function. The slice product of $f$ and $g$ is the slice function on $\Omega_D$ defined by $f{.}g:= \mathcal{I}(FG).$
\end{defi}
In general, $f.g\neq g.f.$ If the components $F_1$, $F_2$ of the first stem function $F$ are real-valued, then $f = \mathcal{I}(F)$ is termed as slice preserving. For the slice preserving function $f$ and slice function $g$, the slice product $f.g$ coincides with $fg$.
\begin{defi}
	For $f=\mathcal{I}(F)\in \mathcal{S}(\Omega_D)$ with $F=F_1+iF_2,$ define the slice conjugate of $f$ as
	\[ f^c= \mathcal{I}\left(\overline{F_1}+i\overline{F_2}\right),\] and the normal function (or symmetrization) of $f$ as 
	\[ N(f)=f.f^c=f^c.f,\] which is slice preserving on $\Omega_D$.
\end{defi}
Let $\mathcal{Z}_f$ denote the zero set of $f$ on $\Omega_D.$ 
\begin{defi}
	Let $f\in \mathcal{S}(\Omega_D).$ If $f$ does not vanish identically, then its slice reciprocal is defined as 
	\[f^{-\bullet}(x):=N(f)(x)^{-1}.f^c(x)=N(f)(x)^{-1}f^(x) \] which is a slice function on $\Omega_D\setminus {\mathcal{Z}_{N(f)}}.$
\end{defi}
More recently, Ghiloni \emph{et. al.} \cite{Ghiloni-Perotti-AdvM-2011} found a new and nice relation between the values of reciprocals $f^{-\bullet}(x)$ and $f^{-1}(x)$ for slice functions $f\in \mathcal{S}(\Omega_D)$ as \[f^{-\bullet}(x)=f\left( T_f(x)\right)^{-1}, \] where $T_f$ is a bijective self-map of $\Omega_D\setminus E$ $E=\{ a+\beta I:I\in \mathcal{S}, z=\alpha+\beta i\in D\;\mbox{for}\; F_2(z)=0\}$ given by 
\[ T_f(x)=\left( f^c(x)^{-1}((xf^c(x)) F_2(z))\right)F_2(z)^{-1},\] which reduces to the known result $T_f(x)= f^c(x)^{-1}xf^c(x)$ for the associative algebra of quaternions.
\section{\bf Proof of main results}\label{Sec-4}
	\begin{proof}[\bf Proof of Theorem \ref{Thm-1.16}]
	For $f(x)=\sum_{k=0}^{\infty}x^ka_k\in \mathcal{SRB}(\mathbb{B})$ and $0<\beta\neq1,$ First we define	
	\begin{align*}
		T^{*,f}_\beta(|x|):=\sum_{k=0}^{\infty}|x|^k\left(\dfrac{1}{k+1}\sum_{n=0}^{k}\dfrac{\Gamma(k-n+\beta)}{\Gamma(k-n+1)\Gamma(\beta)}|a_n|\right).
	\end{align*}
	Let $p=|a_0|<1$ and $|x|=r.$ Under the condition of Theorem \ref{Thm-1.16}, it follows that (see \cite[Theorem 4.1]{Rocchetta-Gentili-Sarfatti-MN-2012})
	\begin{align}\label{Eqq-1.15}
		|a_k|\leq 1-p^2,\;\;\mbox{for}\;\;k\in \mathbb{N}.
	\end{align}
	 This yields
	 \begin{align*}
	 	T^{*,f}_\beta(r)\leq p\sum_{k=0}^{\infty}r^k\left(\dfrac{1}{k+1}\dfrac{\Gamma(k+\beta)}{\Gamma(k+1)\Gamma(\beta)}\right)+(1-p^2)\sum_{k=1}^{\infty}r^k\left(\dfrac{1}{k+1}\sum_{n=1}^{k}\dfrac{\Gamma(k-n+\beta)}{\Gamma(k-n+1)\Gamma(\beta)}\right).
	 \end{align*}
	 The above inequality is equivalent to
	 \begin{align*}
	 		T^{*,f}_\beta(r)&\leq \dfrac{p}{r}\int_{0}^{r}\dfrac{dt}{(1-t)^\beta}+\dfrac{(1-p^2)}{r}\int_{0}^{r}\dfrac{tdt}{(1-t)^{\beta+1}}\\&=\dfrac{(p^2+p-1)}{r}\int_{0}^{r}\dfrac{dt}{(1-t)^\beta}+\dfrac{(1-p^2)}{r}\int_{0}^{r}\dfrac{tdt}{(1-t)^{\beta+1}}\\&\leq \dfrac{1}{r}\left(\dfrac{(p^2+p-1)(1-(1-r)^{1-\beta})}{1-\beta}+\dfrac{(1-p^2)((1-r)^{-\beta}-1)}{\beta}  \right):=\Psi(p).
	 \end{align*}
	 Differentiation of the function $\Psi$ with respect to $p$ gives us
	 \begin{align*}
	 	\begin{cases}
	 		\Psi^\prime(p)=\dfrac{1}{r}\left( \dfrac{(2p+1)\left(1-(1-r)^{1-\beta}\right)}{1-\beta}-\dfrac{2p\left((1-r)^{-\beta}-1 \right)}{\beta} \right),\\\vspace{1.2mm}
	 		\Psi^{\prime \prime}(p)= \dfrac{1}{r}\left( \dfrac{2\left(1-(1-r)^{1-\beta}\right)}{1-\beta}-\dfrac{2\left((1-r)^{-\beta}-1 \right)}{\beta} \right).
	 	\end{cases}
	 \end{align*}
	 It is easy to see that $\Psi^{\prime \prime}(p)\leq 0$ for $p\in [0,1)$ and $r\in [0,1).$ This provides that $\Psi^\prime(p)\geq \Psi^\prime(1)$ and 
	 \[ \Psi^\prime(1)=\dfrac{1}{r}\left( \dfrac{3(1-(1-r)^{1-\beta})}{1-\beta}-\dfrac{2\left( (1-r)^{-\beta}-1\right)}{\beta}\right)\geq 0\] holds for $r\leq R_\beta,$ where $R_\beta$ is the positive root of the equation \eqref{Eq-1.13}. It follows that $\Psi$ is an increasing function of $p$ for $r\leq R_\beta.$ It implies that 
	 \[ \Psi(p)\leq \Psi(1)=\dfrac{1}{r}\left(\dfrac{1-(1-r)^{1-\beta}}{1-\beta}\right),\] for $r\leq R_\beta.$ This is the desired inequality.\vspace{1.2mm}
	 
	 Finally, to show the sharpness, given $p\in [0,1),$ we consider the slice regular function
	 \begin{align}\label{Eqq-1.1}
	 	f_p(x)=(1-xp)^{-\bullet }.(p-x)u=p-(1-p^2)u\sum_{k=1}^{\infty}x^kp^{k-1},\;\;x\in \mathbb{B}
	 \end{align}
	 for some $u\in \partial\mathbb{B}.$ For $f_p$ and $|x|=r,$ we obtain 
	  \begin{align}\label{Eq-1.16}
	 	T^{*,f_p}_\beta(|x|)&= \dfrac{p}{r}\left(\dfrac{(1-(1-r)^{1-\beta})}{1-\beta}\right)+(1-p^2)\sum_{k=1}^{\infty}r^k\left(\dfrac{1}{k+1}\sum_{n=1}^{k}\dfrac{\Gamma(k-n+\beta)}{\Gamma(k-n+1)\Gamma(\beta)}p^{k-1}\right)\\&\nonumber=\dfrac{p}{r}\left(\dfrac{(1-(1-r)^{1-\beta})}{1-\beta}\right)+\dfrac{(1-p^2)}{r}\int_{0}^{r}\dfrac{tdt}{(1-pt)(1-t)^\beta}\\&\nonumber=\dfrac{p}{r}\left(\dfrac{(1-(1-r)^{1-\beta})}{1-\beta}\right)\\&\quad\nonumber-\dfrac{(1-p)}{r}\left(	\dfrac{3(1-(1-r)^{1-\beta})}{1-\beta}-\dfrac{2\left((1-r)^{-\beta}-1\right)}{\beta} \right)+N_p(r),
	 \end{align}
	 where 
	 \begin{align*}
	 	N_p(r):=\dfrac{2(1-p)}{r}\left(	\dfrac{(1-(1-r)^{1-\beta})}{1-\beta}-\dfrac{\left((1-r)^{-\beta}-1\right)}{\beta} \right)+\dfrac{(1-p^2)}{r}\int_{0}^{r}\dfrac{tdt}{(1-pt)(1-t)^\beta}.
	 \end{align*}
	 Using the series representation of $N_p(r)$, we have
	 \begin{align*}
	 		N_p(r)=\sum_{k=0}^{\infty}\dfrac{1}{k+1}\big[&-\dfrac{(1-p)^2}{p}\dfrac{\Gamma(k+\beta)}{\Gamma(k+1)\Gamma(\beta)}-2(1-p)\dfrac{\Gamma(k+\beta+1)}{\Gamma(k+1)\Gamma(\beta+1)}\\&+\dfrac{(1-p^2)}{p}\sum_{m=0}^{k}\dfrac{\Gamma(k-m+\beta)}{\Gamma(k-m+1)\Gamma(\beta)}\big]r^k.
	 \end{align*}
	 By using the identity 
	 \[\sum_{m=0}^{k}\dfrac{\Gamma(k-m+\beta)}{\Gamma(k-m+1)\Gamma(\beta)}=\dfrac{\Gamma(k+\beta+1)}{\Gamma(k+1)\Gamma(\beta+1)} ,\]
	  we can get that $	N_p(r)=O\left((1-p)^2\right),$ as $p$ tends to $1.$ Further, a simple
	 computation shows that for $r>R_\beta$, the quantity
	 \[ \dfrac{3(1-(1-r)^{1-\beta})}{1-\beta}-\dfrac{2\left((1-r)^{-\beta}-1\right)}{\beta}<0.\]
	 
	 After considering these observations in \eqref{Eq-1.16}, we conclude that  $R_\beta$  cannot be improved. This completes the proof.
	\end{proof}
	\begin{proof}[\bf Proof of Theorem \ref{Thm-1.17}]
	Let $f(x)=\sum_{k=m}^{\infty}x^ka_k.\in \mathcal{SRB}(\mathbb{B})$ and 
	\[L^{*,f}_\gamma(|x|):= \sum_{k=m}^{\infty}\dfrac{1}{k+\gamma}|x^ka_k|.\]
	We can write $f(x)=x^mh(x),$ where $h(x):=\sum_{k=m}^{\infty}x^{k-m}a_k.$ Let $p:=|a_m|<1.$ Under the condition of Theorem \ref{Thm-1.17}, it follows that (see \cite[Theorem 4.1]{Rocchetta-Gentili-Sarfatti-MN-2012})
	\begin{align}\label{Eqq-1.155}
		|a_k|\leq 1-p^2,\;\;\mbox{for}\;\;k\geq m+1.
	\end{align}
	In view of \eqref{Eqq-1.155}, we obtain the following inequality
	\[ L^{*,f}_\gamma(|x|)\leq \dfrac{p}{m+\gamma}|x|^m+(1-p^2)\sum_{k=m+1}^{\infty}\dfrac{1}{k+\gamma}|x|^k:=\Phi(p).\]
	We see that 
	\[ \begin{cases}
		\Phi^\prime(p)=\dfrac{1}{m+\gamma}|x|^m-2p\sum_{k=m+1}^{\infty}\dfrac{1}{k+\gamma}|x|^k\vspace{1.2mm}\\
		\Phi^{\prime\prime}(p)=-2\sum_{k=m+1}^{\infty}\dfrac{1}{k+\gamma}|x|^k<0.
	\end{cases}\]
	Thus, 
	\[ 	\Phi^\prime(p)\geq 	\Phi^\prime(1)=\dfrac{1}{m+\gamma}|x|^m-2\sum_{k=m+1}^{\infty}\dfrac{1}{k+\gamma}|x|^k,\]
	for $|x|\leq R_\gamma,$ where $R_\gamma$ is the positive root of the equation \eqref{EEq-1.15}. Hence, $\Phi$ is an increasing function of $p$ for $|x| \leq  R_\gamma$. This gives
		\[ \sum_{k=m}^{\infty}\dfrac{1}{k+\gamma}|x^ka_k|\leq \dfrac{1}{m+\gamma}|x|^m\]
	for $|x|\leq R_\gamma$ and a simple observation gives $R_\gamma<1$.\vspace{1.2mm}
	
	To prove $R_\gamma$ to be the best possible, for $p\in [0,1),$ we consider the function
	\[ g_p(x)=x^m(1-px)^{-\bullet}.(x-p)u=-x^mp+(1-p^2)\sum_{k=1}^{\infty}x^{k+m}up^{k-1}\;\;\mbox{for}\;\; x\in \mathbb{B},\]
	for some $u\in \partial\mathbb{B}.$ For the function $g_p,$ we obtain 
\begin{align*}
	L^{*,g_p}_\gamma(|x|)&= \dfrac{p}{m+\gamma}|x|^m+(1-p^2)\sum_{k=m+1}^{\infty}\dfrac{p^{k-1}}{k+\gamma}|x|^k\\&=\dfrac{p}{m+\gamma}|x|^m-(1-p)\left(\dfrac{|x|^m}{m+\gamma}-2\sum_{k=m+1}^{\infty}\dfrac{|x|^k}{k+\gamma} \right)+M_p(|x|),
\end{align*}
where 
\[M_p(|x|):= 2(p-1)\sum_{k=m+1}^{\infty}\dfrac{|x|^k}{k+\gamma}+(1-p^2)\sum_{k=m+1}^{\infty}\dfrac{p^{k-1}}{k+\gamma}|x|^k. \]
Letting $p\to 1^-,$ we obtain 
\[ M_p(|x|)=\sum_{k=m+1}^{\infty}\dfrac{2(p-1)+(1-p^2)p^{k-1}}{k+\gamma}|x|^k=O\left((1-p)^2\right).\]
Therefore, 
\[ 	L^{*,g_p}_\gamma(|x|)> \dfrac{p}{m+\gamma}|x|^m\] for $|x|>R_\gamma,$ which give that $R_\gamma$ cannot be improved
and the proof is complete.
	\end{proof}
	\begin{proof}[\bf Proof of Theorem \ref{Thm-1.19}]
		 We see that  
		\begin{align*}
			\mathcal{F}^*_f(|x|)&=\sum_{n=0}^{\infty}|x|^n\left(\sum_{k=0}^{n} \bigg|a_ke^{-\frac{2\pi Jnk}{(n+1)}}\bigg|\right)=\sum_{n=0}^{\infty}|x|^n\left(\sum_{k=0}^{n} |a_k|\right).
		\end{align*}
		Let $|a_0|=p.$ In view of \eqref{Eqq-1.155}, we obtain 
		\begin{align*}
			\mathcal{F}^*_f(|x|)\leq p\sum_{n=0}^{\infty}|x|^n+(1-p^2)\sum_{n=1}^{\infty}n|x|^n=\dfrac{p}{1-|x|}+\dfrac{(1-p^2)|x|}{(1-|x|)^2}:=G(p).
		\end{align*}
		Also for $|x|<1,$  we see that
		\[ \begin{cases}
			G^\prime(p)=\dfrac{1}{1-|x|}-\dfrac{2p|x|}{(1-|x|)^2},\vspace{3mm}\\
				G^{\prime\prime}(p)=-\dfrac{2|x|}{(1-|x|)^2}<0.
		\end{cases}\]
		This shows that $	G^\prime$ is a decreasing function of $p$ for $|x|\leq 1/3,$ hence we obtain
		\[ 	G^\prime(p)\geq	G^\prime(1)=\dfrac{(1-3|x|)}{(1-|x|)^2}.\]
		It  can be easily say that $G$ is an increasing function of $p$. Therefore, for $|x|\leq 1/3,$ we see
		\[	G(p)\leq G(1)=\dfrac{1}{1-|x|}.\]This is the desired inequality.\vspace{1.2mm}
		
		We now show that the constant $1/3$ cannot be improved. To
		see this, we consider the function $f_p$ given by \eqref{Eqq-1.1}. For $f_p$, an easy computation leads to
		\begin{align*}
			\mathcal{F}^*_{f_p}(|x|)&= \dfrac{p}{1-|x|}+(1-p^2)\sum_{n=1}^{\infty}|x|^n\left(\sum_{k=1}^{\infty}p^{k-1}\right)\\&=\dfrac{p}{1-|x|}+(1+p)\sum_{n=1}^{\infty}(1-p^n)|x|^n\\&=\dfrac{p}{1-|x|}+(1+p)\left( \dfrac{|x|}{1-|x|}-\dfrac{p|x|}{1-p|x|}\right)\\&=
			\dfrac{1}{1-|x|}-(1-p)\dfrac{1-3|x|}{1-|x|}+H_p(|x|),
		\end{align*}
		where \[H_p(|x|):=\dfrac{(2p-1)2|x|}{1-|x|}-\dfrac{(1+p)p|x|}{1-p|x|}. \]
		Note that $H_p(|x|)\to 0$ as $p\to 1^-$ and $(1-3|x|)/(1-|x|)<0$ when $|x|>1/3.$ This completes the proof.
	\end{proof}
	\begin{proof}[\bf Proof of Theorem \ref{Thm-1.20}]
The proof is similar to Theorem \ref{Thm-1.19}. Let $|a_0|=p.$	In view of \eqref{Eqq-1.155}, we see that 
	\begin{align*}
		\mathcal{L}^*_f(|x|)&\leq p\sum_{n=0}^{\infty}\dfrac{|x|^n}{n+1}+(1-p^2)\sum_{n=1}^{\infty}|x|^n\left(\sum_{k=1}^{n}\dfrac{|p_k|}{(n+1)^{k+1}}\right)\\&=-\dfrac{p}{|x|}\log(1-|x|)+(1-p^2)\sum_{n=1}^{\infty}|x|^n\left( \dfrac{1}{n(1+n)}-\dfrac{1}{n(1+n)^{n+1}}\right)\\&< -\dfrac{p}{|x|}\log(1-|x|)+(1-p^2)\left( \left(\dfrac{1}{|x|}-1\right)\log(1-|x|)+1\right)\\&=\dfrac{1-|x|-p^2+p^2|x|-p}{|x|}\log(1-|x|)+(1-p^2):=Q(p).
	\end{align*}
	We obtain that 
	\[ \begin{cases}
		Q^\prime(p)=\dfrac{-2p+2p|x|-1}{|x|}\log(1-|x|)-2p,\vspace{2mm}\\
		Q^{\prime\prime}(p)= -\dfrac{2(1-|x|)}{|x|}\log(1-|x|)-2<0.
	\end{cases}\]
	This shows that $Q^\prime$ is a decreasing function, which gives
	\[ Q^\prime(p)\geq Q^\prime(1)=\dfrac{2|x|-3}{|x|}\log(1-|x|)-2>0\] for all $|x|<1$. Thus, we have
\begin{align*}
	Q(p)\leq Q(1)=\dfrac{1}{|x|}\log \frac{1}{(1-|x|)}. 
\end{align*}
	This completes the proof.
	\end{proof}
	\begin{proof}[\bf Proof of Corollary \ref{Cor-2.3}]
	We now reformulate the right-hand side of inequality \eqref{Eq-3.3} as follows:
		\begin{align*}
			\mathcal{L}^{*}_f(r)&\leq \frac{1}{r}\ln\left(\frac{1}{1-r}\right)+\frac{(1-p)}{r}\ln (1-r)+ 2(1-p)\left(\frac{1}{r}\ln\left(\frac{1}{1-r}\right)-\frac{{\rm Li}_2(r)}{r}\right)\\&\leq\frac{1}{r}\ln\left(\frac{1}{1-r}\right)+(1-p)\bigg[\frac{1}{r}\ln(1-r)+2\left(\frac{1}{r}\ln\left(\frac{1}{1-r}\right)-\frac{{\rm Li}_2(r)}{r}\right)\bigg]\\&\leq \frac{1}{r}\ln\left(\frac{1}{1-r}\right)+(1-p)\left(\frac{1}{r}\ln(1-r)+2\left(\frac{1}{r}\ln\left(\frac{1}{1-r}\right)-\frac{{\rm Li}_2(r)}{r}\right)\right).
		\end{align*}
		Thus, the desired inequality
		\begin{align}\label{Eq-3.5}
			\mathcal{L}^{*}_f(r)\leq \frac{1}{r}\ln\left(\frac{1}{1-r}\right)
		\end{align}
		can be obtained if $\mathcal{P}(r)\leq 0$ for $r\leq r^*<1,$ where $\mathcal{P}(r)$ is as in the statement.  Applying fundamental theorems of calculus, we obtain $r^*$ as a root of the equation $\mathcal{P}(r) = 0$ in the interval $(0,1)$. Consequently, $r^*$ serves as the Bohr radius for the discrete Laplace transform $\mathcal{L}^*_f$.\vspace{1.2mm}
		
		Moreover, it can be shown that the inequality \eqref{Eq-3.5} is sharp and this can be shown using the function $f_p$ as given by \eqref{Eqq-1.1}. For this $f_p$, a tedious computation gives that
		\begin{align*}
			\mathcal{L}^{*}_{f_{p}}(r)&=|a_0|\sum_{n=0}^{\infty}\frac{|x|^n}{n+1}+\sum_{n=1}^{\infty}|x|^n\left(\sum_{k=1}^{n}\frac{|a_k|}{(n+1)^{k+1}}\right)\\&=\frac{p}{r}\ln\left(\frac{1}{1-r}\right)+\frac{1-p^2}{p}\sum_{n=1}^{\infty}\left(\sum_{k=1}^{n}\frac{p^k}{(n+1)^{k+1}}\right)r^k\\&=\frac{1}{r}\ln\left(\frac{1}{1-r}\right)+(1-p)\bigg[-\frac{1}{r}\ln\left(\frac{1}{1-r}\right)\\&\quad+\frac{1+p}{p}\sum_{n=1}^{\infty}\left(\sum_{k=1}^{n}\frac{p^k}{(n+1)^{k+1}}\right)r^k\bigg].
		\end{align*}
		Taking $p\to 1^{-}$ in the last expression, it is easy to see that 
		\begin{align*}
			\mathcal{L}^{*}_{f_{p}}(r)=\frac{1}{r}\ln\left(\frac{1}{1-r}\right)
		\end{align*}
		which shows that the bound is sharp. This completes the proof.
	\end{proof}		
	\noindent{\bf Acknowledgment:} The research of the first author is supported by UGC-JRF (Ref. No. 201610135853), New Delhi, Govt. of India and second author is supported by SERB File No. SUR/2022/002244, Govt. of India.\vspace{2mm}

	\noindent {\bf Funding:} Not Applicable.\vspace{1.5mm}
	
	\noindent\textbf{Conflict of interest:} The authors declare that there is no conflict  of interest regarding the publication of this paper.\vspace{1.5mm}
	
	\noindent\textbf{Data availability statement:}  Data sharing not applicable to this article as no datasets were generated or analysed during the current study.\vspace{1.5mm}
	
	\noindent{\bf Code availability:} Not Applicable.\vspace{1.5mm}
	
	\noindent {\bf Authors' contributions:} Both the authors have equal contributions.
  
\end{document}